\documentclass[11pt,a4paper]{article}

\usepackage[margin=1.2in]{geometry}
\usepackage[utf8]{inputenc}
\usepackage[english]{babel}
\usepackage{booktabs}
\usepackage{amsmath}
\usepackage{natbib}
\usepackage{algorithm, algorithmicx, algpseudocode}
\usepackage{graphicx}
\usepackage{pgfplots}
\usepackage{pgfplotstable}
\pgfplotsset{width=10cm,compat=1.9}
\usetikzlibrary{patterns}
\pgfplotsset{compat = 1.13,
        my ybar legend/.style={
            legend image code/.code={
                \draw [##1] (0cm,-0.6ex) rectangle +(2em,1.5ex);
            },
        },
}

\newenvironment{tablenotes}{\list{}{\setlength{\labelsep}{0pt}%
\setlength{\labelwidth}{0pt}%
\setlength{\leftmargin}{0pt}%
\setlength{\rightmargin}{0pt}%
\setlength{\topsep}{2pt}%
\setlength{\itemsep}{0pt}%
\setlength{\partopsep}{0pt}%
\setlength{\listparindent}{0em}%
\setlength{\parsep}{0pt}}%
\item\relax%
}{\endlist}%

\usepackage[scientific-notation=true]{siunitx}

\usepackage{caption}
\usepackage{subcaption}

\title{A new deflation criterion for the QZ algorithm\footnotemark[1]}

\author{Thijs Steel, Raf Vandebril, Julien Langou}

\begin{document}

\maketitle

\abstract{The QZ algorithm computes the generalized Schur form of a matrix pencil. It is an iterative algorithm and, at some point, it must decide when to deflate, that is when a generalized eigenvalue has converged and to move on to another one. Choosing a deflation criterion that makes this decision is nontrivial. If it is too strict, the algorithm might waste iterations on already converged eigenvalues. If it is not strict enough, the computed eigenvalues might not have full accuracy. Additionally, the criterion should not be computationally expensive to evaluate. There are two commonly used criteria:
the \emph{elementwise} criterion and the \emph{normwise} criterion. This paper introduces a new deflation criterion based on the size of and the gap between the eigenvalues. We call this new deflation criterion the \emph{strict} criterion. This new criterion for QZ is analogous to the criterion derived by Ahues and Tisseur~\cite{ahues1997new} for the QR algorithm. Theoretical arguments and numerical experiments suggest that 
the strict criterion
outperforms the normwise and elementwise criteria in terms of accuracy. We also provide
an example where the accuracy of the generalized eigenvalues
using the elementwise or the normwise criteria
is less than 2 digits whereas the strict criterion 
leads to generalized eigenvalues
which are almost accurate to the working precision.
Additionally, this paper evaluates some commonly used criteria for infinite eigenvalues.}

\footnotetext[1]{\textbf{Abbreviations:} IEEE, Institute of Electrical and Electronics Engineers; LAPACK, Linear Algebra Package}
\footnotetext[2]{\textbf{Published version:} https://doi.org/10.1002/nla.2524}

\section{Introduction}

To solve dense generalized eigenvalue problems, the most popular algorithm is the QZ algorithm~\cite{moler1973algorithm}. It consists of a direct reduction to Hessenberg-triangular form, an iterative reduction to generalized Schur form, and finally a direct step to calculate the eigenvectors. This text focuses on the iterative reduction to generalized Schur form.

Consider a pencil $(H,T)$ in Hessenberg-triangular form:
\begin{equation}
   (H,T) = 
   \left( \begin{bmatrix}
      h_{11} & \dots  & \dots     & h_{1,n} \\
      h_{21} & \ddots &           & \vdots  \\
             & \ddots & \ddots    & \vdots  \\
      0      &        & h_{n,n-1} & h_{nn}  \\
   \end{bmatrix}, 
   \begin{bmatrix}
      t_{11} & \dots  & \dots  & t_{1,n} \\
             & \ddots &        & \vdots  \\
             &        & \ddots & \vdots  \\
      0      &        &        & t_{nn}  \\
   \end{bmatrix} \right).
   \label{eq:pencil}
\end{equation}
On a high level, the iterative reduction of this pencil to Schur form consists of the following steps:
\begin{enumerate}
   \item Select shifts. Often, Wilkinson shifts are used: the eigenvalues of the $2\times 2$ subpencil in the bottom right.
   \item Chase a bulge through the pencil using these shifts.
   \item If a subdiagonal entry in $H$ is small enough, it is set to zero and the eigenvalue problem is split into separate subproblems.
   \item If a diagonal entry in $T$ is small enough, it is set to zero, and the pencil is transformed to move this entry to the top or bottom of the pencil. Once at the edge, a single infinite eigenvalue can be split off from the pencil.
   \item Repeat steps 1 through 4 until every subproblem is of size $1\times 1$ (or in the case of complex eigenvalues in real pencils: $2\times 2$).
\end{enumerate}
More detailed explanations of the QZ algorithm can be found in \cite{moler1973algorithm, kressnerbook}. This paper focuses on steps 3 and 4, and more precisely on the deflation criterion to determine whether the elements are small enough that they can be set to zero. If this criterion is too strict, we may need to do extra iterations when the computed eigenvalues were already accurate enough. If the criterion is not strict enough, we may get inaccurate eigenvalues.

Moler and Stewart~\cite{moler1973algorithm}, in the paper where they introduced the QZ algorithm for the first time, consider a subdiagonal entry small enough if $|h_{i+1,i}| \le u\|H\|_F$ and consider a diagonal entry small enough if $|t_{i,i}| \le u\|T\|_F$. Here and throughout the rest of this text, $u$ will refer to the machine precision\footnote{There are different definitions of the machine precision in literature. We define $u$ as the distance between 1 and the smaller floating point number larger than 1.}. We will refer to this as the \emph{normwise criterion}. If we want a backward stable decomposition, we cannot set any entries that do not satisfy this criterion to zero.

In the QR algorithm, it is known that a similar normwise criterion leads to suboptimal results for graded matrices \cite{stewart1990eigenvalues}. Accuracy for such matrices is improved by using a criterion we will refer to as \emph{elementwise}. This criterion was also adapted to the QZ algorithm \cite{adlerborn2006parallel, kaagstrom2007multishift, steel2021multishift}, although we are not aware of detailed research supporting it. With this criterion, a subdiagonal entry is considered small enough if $|h_{i+1,i}| \le u(|h_{i,i}| + |h_{i+1,i+1}|)$ and a diagonal entry is considered small enough if $|t_{i,i}| \le u(|t_{i-1,i}| + |t_{i,i+1}|)$. Note that $|h_{i,i}| + |h_{i+1,i+1}| \le \sqrt{2}\sqrt{|h_{i,i}|^2 + |h_{i+1,i+1}|^2} \le \sqrt{2}\|H\|_F$ so an entry that satisfies the elementwise criterion will also satisfy the normwise criterion up to a constant factor of $\sqrt{2}$. Additionally, $|h_{i,i}| + |h_{i+1,i+1}|$ is usually much smaller than $\|H\|_F$ because only two entries are considered. In the rest of this text, we will ignore this constant factor of $\sqrt{2}$ and say that the elementwise criterion is stricter than the normwise criterion.

The linear algebra library LAPACK \cite{anderson1999lapack} originally used the normwise criterion 
for finite and infinite eigenvalues. It switched to the elementwise criterion in version 3.9.1
for finite and infinite eigenvalues. And, as of version 3.11.0, the latest version at the time of writing, it switched to the elementwise criterion for finite eigenvalues and the normwise criterion for infinite eigenvalues.

The deflation criterion in the QR algorithm has received attention in the work of Ahues and Tisseur~\cite{ahues1997new}. Their paper provides a theoretical basis to support a different deflation criterion than normwise and elementwise. They approximate the sensitivity of the eigenvalues w.r.t. the subdiagonal elements and use this approximation to argue that a subdiagonal element is small enough if it satisfies:
\begin{equation*}
   |h_{i,i-1}| |h_{i-1,i}| \le u | h_{i,i} - h_{i-1,i-1} | \text{ and } |h_{i+1,i}| \le u(|h_{i,i}| + |h_{i+1,i+1}|).
\end{equation*}
So they derived and argue for a new 
deflation criterion for the QR algorithm, which has also found its way into LAPACK.
This paper aims to extend their work to the QZ algorithm. 
We will call our new deflation criterion for QZ the 
\emph{strict} criterion.

\section{Deflation of finite eigenvalues}

In this section, we assume that the pencil only has finite eigenvalues.

\subsection{Derivation}

To derive our new criterion, we will focus our attention on a $2\times 2$ subpencil of Equation~\eqref{eq:pencil}:

\begin{equation}
   (H_i(\epsilon), T_i) = \left(\begin{bmatrix}
      h_{i-1,i-1}      & h_{i-1,i} \\
      \epsilon & h_{i,i}   \\
   \end{bmatrix}, 
   \begin{bmatrix}
      t_{i-1,i-1} & t_{i-1,i} \\
      0           & t_{i,i} \\
   \end{bmatrix} \right).
\end{equation}
Let $\lambda(\epsilon)$ be the eigenvalue of the subpencil closest to $h_{i,i}/t_{i,i}$. It is important to note that $\lambda(\epsilon)$ is not necessarily an eigenvalue of the full pencil. In the rest of this analysis, we will assume that the influence of entries outside this subpencil is negligible. At this point we would like to stress that this is not a detailed error analysis. We are designing a heuristic and significant terms will be approximated.

We want to know when $\epsilon$ is small enough to be considered zero. This is the case when
\begin{equation}\label{eq:forwarderrorcheck}
   |\lambda(\epsilon) - \lambda(0)| \le u |\lambda(0)|,
\end{equation}
where $u$ is the machine precision. We do not know $\lambda(\epsilon)$, so we must approximate.
If the function $\lambda(\epsilon)$ is analytic, then a first order expansion gives the bound:
\begin{equation}\label{eq:forwarderrorbound}
   |\lambda(\epsilon) - \lambda(0)| = |\lambda'(0)\epsilon| + O(|\epsilon|^2).
\end{equation}
To calculate $|\lambda'(0)|$, we first calculate the characteristic polynomial:
\begin{equation*}
   \begin{split}
      p(\epsilon, x) =& \det( H_i(\epsilon) - x T_i ) \\
      =& (h_{i-1,i-1} - xt_{i-1,i-1})(h_{i,i}-xt_{i,i}) - (h_{i-1,i} - xt_{i-1,i})\epsilon. \\
   \end{split}.
\end{equation*}
The characteristic polynomial evaluated in one of the eigenvalues is always zero, so $\frac{d}{d\epsilon}p(\epsilon,\lambda(\epsilon)) = 0$ and we can extract an expression for $\lambda'(\epsilon)$:
\begin{equation}\label{eq:dpexpansion}
   \begin{split}
      \frac{d}{d\epsilon}p(\epsilon,\lambda(\epsilon)) &= 0\\
      \Rightarrow \left. \frac{\partial p}{\partial \epsilon}(\epsilon,x)\right\vert_{x = \lambda(\epsilon)} + \lambda'(\epsilon)\left.\frac{\partial p}{\partial x}(\epsilon,x)\right\vert_{x = \lambda(\epsilon)} &= 0\\
      \Rightarrow \lambda'(\epsilon) &= \dfrac{\left. \frac{\partial p}{\partial \epsilon}(\epsilon,x)\right\vert_{x = \lambda(\epsilon)}}{\left.\frac{\partial p}{\partial x}(\epsilon,x)\right\vert_{x = \lambda(\epsilon)}}.
   \end{split}
\end{equation}
Next, we extract an expression for $\frac{\partial p}{\partial x}(\epsilon,x)$.
\begin{equation*}
   \begin{split}
      \frac{\partial p}{\partial x}(\epsilon,x) =& -t_{i-1,i-1}(h_{i,i}-xt_{i,i})\\
      & - t_{i,i}(h_{i-1,i-1} - xt_{i-1,i-1})\\
      & + t_{i-1,i}\epsilon.\\
   \end{split}.
\end{equation*}
Luckily, almost all terms in this expression are zero when evaluated for $\epsilon = 0$ and $x = h_{i,i}t_{i,i}^{-1}$.
\begin{equation}\label{eq:dpdx}
   \begin{split}
      \frac{\partial p}{\partial x}(0,h_{i,i}t_{i,i}^{-1}) =& - (h_{i-1,i-1}t_{i,i} - h_{i,i}t_{i-1,i-1}).\\
   \end{split}
\end{equation}
Next, we extract an expression for $\frac{\partial p}{\partial \epsilon}(\epsilon,x)$:
\begin{equation*}
   \begin{split}
      \frac{\partial p}{\partial \epsilon}(\epsilon,x) =& - (h_{i-1,i} - xt_{i-1,i}).\\
   \end{split}
\end{equation*}
Evaluated for $\epsilon = 0$ and $x = h_{i,i}t_{i,i}^{-1}$, this becomes:
\begin{equation}\label{eq:dpdeps}
   \begin{split}
      \frac{\partial p}{\partial \epsilon}(0,h_{i,i}t_{i,i}^{-1}) =& - t_{i,i}^{-1}(h_{i-1,i}t_{i,i} - h_{i,i}t_{i-1,i})\\
   \end{split}.
\end{equation}
Combining Equations \eqref{eq:dpexpansion}, \eqref{eq:dpdx} and \eqref{eq:dpdeps} leads to:
\begin{equation}
   \lambda'(0) = \frac{h_{i-1,i}t_{i,i} - h_{i,i}t_{i-1,i}}{t_{i,i}(h_{i-1,i-1}t_{i,i} - h_{i,i}t_{i-1,i-1})}.
   \label{eq:lambdaderivative}
\end{equation}
We can follow the same derivation for $\bar{\lambda}$, the eigenvalue closest to $\frac{h_{i-1,i-1}}{t_{i-1,i-1}}$ to get the similar expression:
\begin{equation}
   \bar{\lambda}'(0) = \frac{h_{i-1,i-1}t_{i-1,i} - h_{i-1,i}t_{i-1,i-1}}{t_{i-1,i-1}(h_{i-1,i-1}t_{i,i} - h_{i,i}t_{i-1,i-1})}.
   \label{eq:lambdaderivative2}
\end{equation}
To design a deflation criterion using these equations, we will combine Equation \eqref{eq:lambdaderivative}, Equation \eqref{eq:forwarderrorcheck} and \eqref{eq:forwarderrorbound}. This leads to:
\begin{equation}
   |h_{i-1,i}t_{i,i} - h_{i,i}t_{i-1,i}||h_{i,i-1}| \le u |h_{i,i}||h_{i-1,i-1}t_{i,i} - h_{i,i}t_{i-1,i-1}|.
   \label{eq:finitedeflation}
\end{equation}
Finally, imagine a pencil where $h_{i-1,i} = t_{i-1,i} = 0$. The eigenvalues of the subpencil will now be independent of $\epsilon$, so the deflation criterion will be satisfied, but setting $\epsilon$ to zero could lead to a large backward error. To guarantee backward stability the criterion will also check that:
\begin{equation*}
   |h_{i,i-1}| \le u (|h_{i,i}| + |h_{i-1,i-1}|).
\end{equation*}

\subsection{Experiments}
In these experiments, we will evaluate the accuracy of the QZ algorithm while using the following criteria:
\begin{itemize}
   \item normwise: $|h_{i,i-1}| \le u \|H\|_F$
   \item elementwise: $|h_{i,i-1}| \le u (|h_{i-1,i-1}| + |h_{i,i}|)$
   \item strict (our new criterion): $|h_{i,i-1}| \le u (|h_{i-1,i-1}| + |h_{i,i}|)$ and $|h_{i-1,i}t_{i,i} - h_{i,i}t_{i-1,i}||h_{i,i-1}| \le u |h_{i,i}||h_{i-1,i-1}t_{i,i} - h_{i,i}t_{i-1,i-1}|$
\end{itemize}
Our implementation consists of a lightly edited version of XHGEQZ, the single/double shift QZ implementation in LAPACK \cite{anderson1999lapack}.

Before describing the different test cases, we will describe a procedure to generate a pencil with predetermined eigenvalues and a predetermined condition number of the eigenvector matrix $\kappa$. This procedure is described in Algorithm \ref{alg:testpencil}. We stress that the condition of the eigenvector matrix is not necessarily related to the condition of the eigenvalues.

\begin{algorithm}[h!]
  \caption{Generation of test pencil}\label{alg:testpencil}
  \begin{algorithmic}
    \State Make $D_1$ and $D_2$ diagonal matrices, setting the diagonal entries to the desired eigenvalue pairs $(\alpha_i,\beta_i)$. (Generalized eigenvalues are often described as pairs, with $\lambda_i = \frac{\alpha_i}{\beta_i})$.
    \State Generate diagonal matrices $V$ and $W$ with diagonal entries logarithmically spaced between 1 and $\kappa^{-1}$ ($\kappa$ determines the condition number of the matrix that contains the eigenvectors).
    \State Multiply $V$ and $W$ from the left and right with random unitary matrices (we choose to generate such matrices by taking QR factorizations of matrices with uniformly random entries).
    \State Generate the final pencil as $(A,B) = V(D_1,D_2)W$.
  \end{algorithmic}
  \end{algorithm}


To show that the new deflation criterion can improve the accuracy of the computed eigenvalues significantly, we propose a small example similar to the one by Ahues and Tisseur \cite{ahues1997new}. Consider the pencil:
\begin{equation}
   (H,T) = \left(
   \begin{bmatrix}
      1 & c & 0 \\
      \eta & (1 + d) & 1 \\
      0 & \eta & (1+ 2d)c^{-1}
   \end{bmatrix},
   \begin{bmatrix}
      1 & 0 & 0 \\
      0 & 1 & 0 \\
      0 & 0 & c^{-1}
   \end{bmatrix}
   \right),
\end{equation}
with $\eta = 1.1 \cdot 10^{-8}$, $c = 1.1 \cdot 10^5$ and $d=10^{-2}$. The eigenvalues of this pencil are approximately: $0.95371503$, $1.0261424$, and $1.0501426$. This pencil shows that a good criterion needs to take both $H$ and $T$ into account. Criteria applied only to $H$ (the normwise criterion, the elementwise criterion or even the criterion by Ahues and Tisseur \cite{ahues1997new}) would consider $h_{3,2}$ to be small enough in IEEE single precision. The returned eigenvalues are then $1$, $1.02$ and $1.04$, resulting in less than 2 digits of accuracy. When using our new criterion, the QZ algorithm performs 2 more iterations, resulting in smaller subdiagonal entries. The returned eigenvalues are: $0.95371503$, $1.0261424$, and $1.0501425$, which are almost accurate to the working precision.

From this example, we can already conclude that there are cases where our new criterion significantly improves the accuracy of the eigenvalues. It is also very unlikely that the other criteria would lead to better results because our new criterion (if we include the elementwise criterion for backward stability) is stricter than the other criteria.

We now test the criteria on matrices of certain classes with randomness in them:

\begin{itemize}
   \item Unitarily diagonalizable pencils. For these pencils, we use Algorithm \ref{alg:testpencil} with the $\alpha_i$ randomly generated, $\beta_i = 1$ and $\kappa = 1$.
   \item Non-unitarily diagonalizable pencils. For these pencils, we use Algorithm \ref{alg:testpencil} with the $\alpha_i$ randomly generated, $\beta_i = 1$ and $\kappa = 1,000$.
   \item Graded matrices. These pencils are generated by first generating two full random matrices $A$ and $B$ and then scaling the rows and columns with values varying between $1$ and $10^{-3}$.
\end{itemize}

\begin{figure}
   \centering
   \begin{subfigure}[b]{0.49\textwidth}
      \resizebox{\textwidth}{!}{%
         \begin{tikzpicture}
            \begin{axis}[
                  ybar,
                  bar width=2mm,
                  area style,
                  xmin = 3.5,
                  xtick={4, 5, 6, 7},
                  xticklabels={$10^{-4}$,$10^{-5}$,$10^{-6}$,$10^{-7}$},
                  xlabel=Relative error,
                  ylabel=Bin count,
                  legend style={
                        at={(0,1)},
                        anchor=north west,
                     },
               ]
               \addplot+[ybar,mark=no,color=black,bar shift=-3mm,fill=black!10,postaction={pattern=dots}] table [x=bins, y=normwise, col sep=comma] {Figures/accuracy_well_conditioned.csv};
               \addplot+[ybar,mark=no,color=black,bar shift=-1mm,fill=black!50,postaction={pattern=north east lines}] table [x=bins, y=elementwise, col sep=comma] {Figures/accuracy_well_conditioned.csv};
               \addplot+[ybar,mark=no,color=black,bar shift=1mm, fill=black!50] table [x=bins, y=strict, col sep=comma] {Figures/accuracy_well_conditioned.csv};
               \legend{normwise,elementwise,strict}
            \end{axis}
         \end{tikzpicture}
      }
   \end{subfigure}
   \hfill
   \begin{subfigure}[b]{0.49\textwidth}
      \resizebox{\textwidth}{!}{%
         \begin{tikzpicture}
            \begin{axis}[
                  ybar,
                  bar width=2mm,
                  area style,
                  xlabel=Number of iterations,
                  ylabel=Bin count
               ]
               \addplot+[ybar,mark=no,color=black,bar shift=-3mm,fill=black!10,postaction={pattern=dots}] table [x=bins, y=normwise, col sep=comma] {Figures/iterations_well_conditioned.csv};
               \addplot+[ybar,mark=no,color=black,bar shift=-1mm,fill=black!50,postaction={pattern=north east lines}] table [x=bins, y=elementwise, col sep=comma] {Figures/iterations_well_conditioned.csv};
               \addplot+[ybar,mark=no,color=black,bar shift=1mm, fill=black!50] table [x=bins, y=strict, col sep=comma] {Figures/iterations_well_conditioned.csv};
            \end{axis}
         \end{tikzpicture}
      }
   \end{subfigure}
   \caption{Histogram of the maximum relative errors (left) and the number of QZ iterations (right) for the eigenvalues of randomly generated pencils with small $\kappa$ for different deflation criteria.}
   \label{fig:expnormal}
\end{figure}
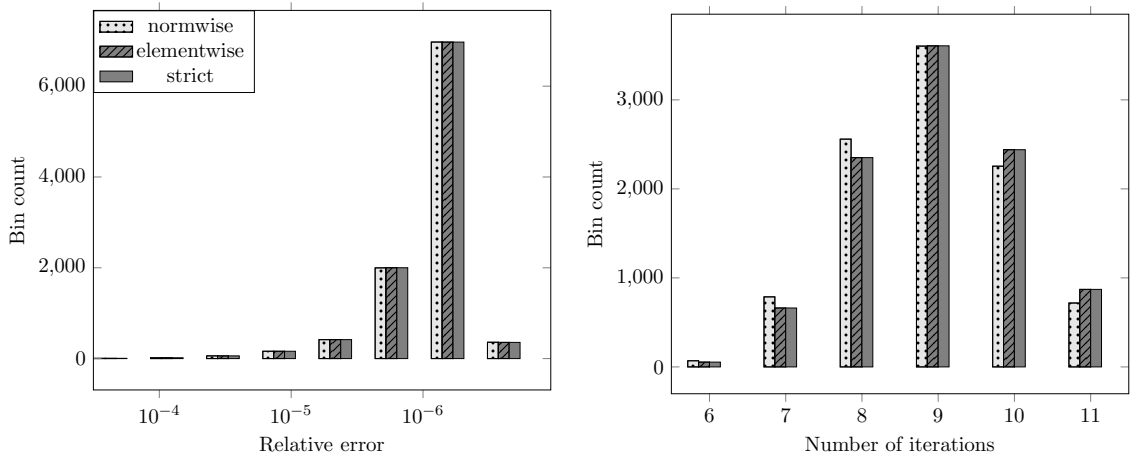

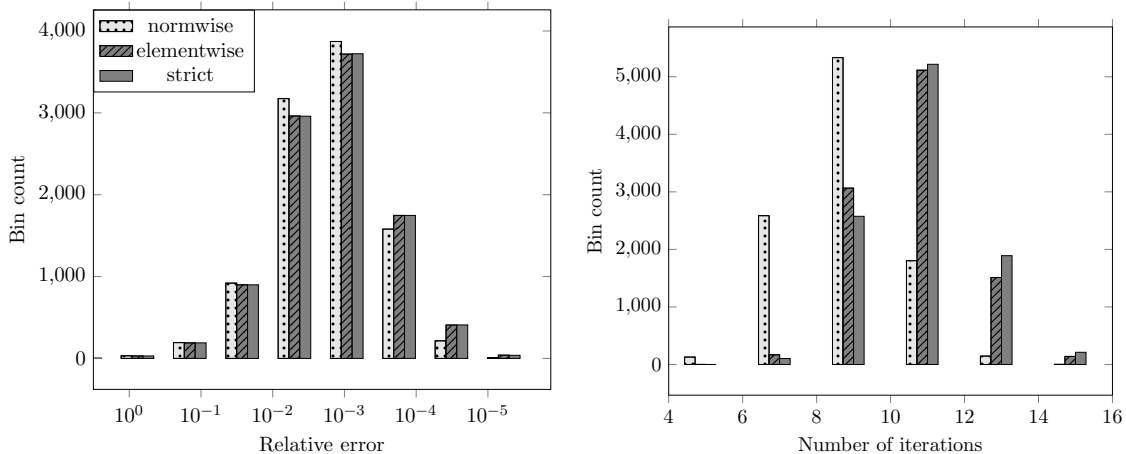
\begin{figure}
   \centering
   \begin{subfigure}[b]{0.49\textwidth}
      \resizebox{\textwidth}{!}{%
         \begin{tikzpicture}
            \begin{axis}[
                  ybar,
                  bar width=2mm,
                  area style,
                  xmin = -0.5,
                  xtick={0, 1, 2, 3, 4, 5, 6, 7},
                  xticklabels={$10^{0}$,$10^{-1}$,$10^{-2}$,$10^{-3}$,$10^{-4}$,$10^{-5}$,$10^{-6}$,$10^{-7}$},
                  xlabel=Relative error,
                  ylabel=Bin count,
                  legend style={
                        at={(0,1)},
                        anchor=north west,
                     },
               ]
               \addplot+[ybar,mark=no,color=black,bar shift=-3mm,fill=black!10,postaction={pattern=dots}] table [x=bins, y=normwise, col sep=comma] {Figures/accuracy_ill_conditioned.csv};
               \addplot+[ybar,mark=no,color=black,bar shift=-1mm,fill=black!50,postaction={pattern=north east lines}] table [x=bins, y=elementwise, col sep=comma] {Figures/accuracy_ill_conditioned.csv};
               \addplot+[ybar,mark=no,color=black,bar shift=1mm, fill=black!50] table [x=bins, y=strict, col sep=comma] {Figures/accuracy_ill_conditioned.csv};
               \legend{normwise,elementwise,strict}
            \end{axis}
         \end{tikzpicture}
      }
   \end{subfigure}
   \hfill
   \begin{subfigure}[b]{0.49\textwidth}
      \resizebox{\textwidth}{!}{%
         \begin{tikzpicture}
            \begin{axis}[
                  ybar,
                  bar width=2mm,
                  area style,
                  xlabel=Number of iterations,
                  ylabel=Bin count
               ]
               \addplot+[ybar,mark=no,color=black,bar shift=-3mm,fill=black!10,postaction={pattern=dots}] table [x=bins, y=normwise, col sep=comma] {Figures/iterations_ill_conditioned.csv};
               \addplot+[ybar,mark=no,color=black,bar shift=-1mm,fill=black!50,postaction={pattern=north east lines}] table [x=bins, y=elementwise, col sep=comma] {Figures/iterations_ill_conditioned.csv};
               \addplot+[ybar,mark=no,color=black,bar shift=1mm, fill=black!50] table [x=bins, y=strict, col sep=comma] {Figures/iterations_ill_conditioned.csv};
            \end{axis}
         \end{tikzpicture}
      }
   \end{subfigure}
   \caption{Histogram of the maximum relative errors (left) and the number of QZ iterations (right) for the eigenvalues of randomly generated non-unitarily diagonalizable pencils for different deflation criteria.}
   \label{fig:expnonnormal}
\end{figure}

\begin{figure}
   \centering
   \begin{subfigure}[b]{0.49\textwidth}
      \resizebox{\textwidth}{!}{%
         \begin{tikzpicture}
            \begin{axis}[
                  ybar,
                  bar width=2mm,
                  area style,
                  xmin=-1.5,
                  xtick={-1,0, 1, 2, 3, 4, 5, 6, 7},
                  xticklabels={$10^{1}$,$10^{0}$,$10^{-1}$,$10^{-2}$,$10^{-3}$,$10^{-4}$,$10^{-5}$,$10^{-6}$,$10^{-7}$},
                  xlabel=Relative error,
                  ylabel=Bin count,
                  legend style={
                        at={(0,1)},
                        anchor=north west,
                     },
               ]
               \addplot+[ybar,mark=no,color=black,bar shift=-3mm,fill=black!10,postaction={pattern=dots}] table [x=bins, y=normwise, col sep=comma] {Figures/accuracy_graded.csv};
               \addplot+[ybar,mark=no,color=black,bar shift=-1mm,fill=black!50,postaction={pattern=north east lines}] table [x=bins, y=elementwise, col sep=comma] {Figures/accuracy_graded.csv};
               \addplot+[ybar,mark=no,color=black,bar shift=1mm, fill=black!50] table [x=bins, y=strict, col sep=comma] {Figures/accuracy_graded.csv};
               \legend{normwise,elementwise,strict}
            \end{axis}
         \end{tikzpicture}
      }
   \end{subfigure}
   \hfill
   \begin{subfigure}[b]{0.49\textwidth}
      \resizebox{\textwidth}{!}{%
         \begin{tikzpicture}
            \begin{axis}[
                  ybar,
                  bar width=2mm,
                  area style,
                  xlabel=Number of iterations,
                  ylabel=Bin count
               ]
               \addplot+[ybar,mark=no,color=black,bar shift=-3mm,fill=black!10,postaction={pattern=dots}] table [x=bins, y=normwise, col sep=comma] {Figures/iterations_graded.csv};
               \addplot+[ybar,mark=no,color=black,bar shift=-1mm,fill=black!50,postaction={pattern=north east lines}] table [x=bins, y=elementwise, col sep=comma] {Figures/iterations_graded.csv};
               \addplot+[ybar,mark=no,color=black,bar shift=1mm, fill=black!50] table [x=bins, y=strict, col sep=comma] {Figures/iterations_graded.csv};
            \end{axis}
         \end{tikzpicture}
      }
   \end{subfigure}
   \caption{Histogram of the maximum relative errors (left) and the number of QZ iterations (right) for the eigenvalues of randomly generated graded pencils for different deflation criteria.}
   \label{fig:expgraded}
\end{figure}
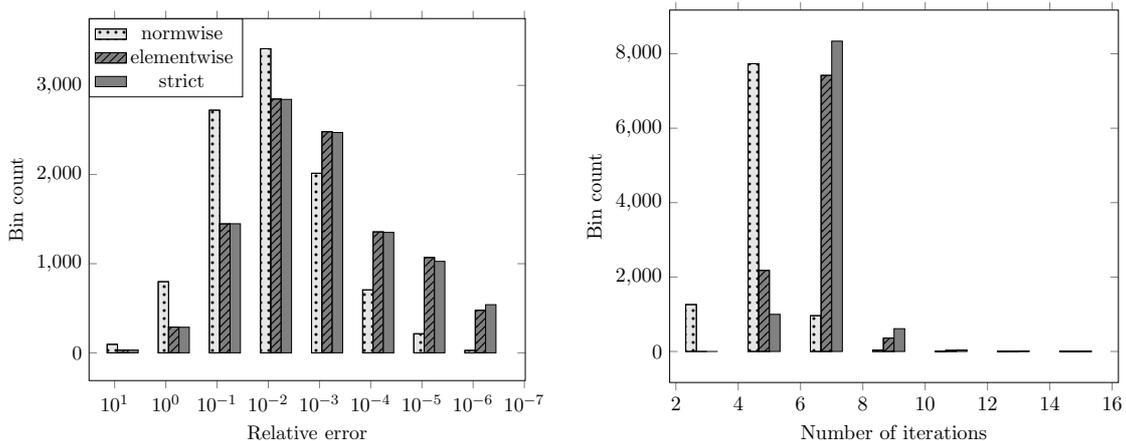

To perform the test, the pencils are reduced to Hessenberg-triangular form, and then the eigenvalues are computed in double precision and compared with the eigenvalues computed in single precision. For each class, 10,000 pencils of size 50 were generated.
Histograms of the maximum relative error of the eigenvalues for the three classes of matrices with the three different types of deflation criteria are shown in Figure \ref{fig:expnormal}, \ref{fig:expnonnormal} and \ref{fig:expgraded}. The averages of the distributions are shown in Table \ref{table: averages}.

\begin{center}
  \begin{table*}[t]%
  \caption{Average number of accurate digits in the eigenvalues for different classes of pencils and different deflation criteria. \label{table: averages}}
  \centering
  \begin{tabular*}{400pt}{@{\extracolsep\fill}lccc@{\extracolsep\fill}}
  \toprule
  & \textbf{normwise}  & \textbf{elementwise}  & \textbf{strict}  \\
  \midrule
  \textbf{Unitarily diagonalizable}    & 6.28 (8.9)  & 6.28 (9.1)   & 6.28 (9.1) \\
  \textbf{Non-unitarily diagonalizable} & 3.20 (9.4) & 3.27 (11.2)   & 3.27 (11.4) \\
  \textbf{graded}            & 2.70 (5.5)  & 3.56 (7.1)    & 3.57 (7.4)  \\
  \bottomrule
  \end{tabular*}
  \begin{tablenotes}
  \item The average number of required iterations is shown between brackets.
  \end{tablenotes}
  \end{table*}
\end{center}

We notice a significant increase in accuracy for graded pencils. If we look at the histogram, there is a slight difference in accuracy when using the strict criterion over the elementwise criterion. However, performing a t-test on the logarithm of the accuracy of the eigenvalues of the non-unitarily diagonalizable pencils reveals that we cannot conclude that the strict criterion results in higher accuracy on average $(p = 0.704)$. No convergence failures occurred.

These experiments lead us to conclude that while there are certainly pencils where the strict criterion leads to a significant increase in accuracy over the elementwise criterion, the average accuracy of the two criteria is very similar. Both the elementwise and strict criterion result in a significant increase in accuracy for graded pencils. We cannot rule out that a class of generalized eigenvalue problems for which there would be an increase in the average accuracy exists. For a general-purpose library like LAPACK, we believe changing to the strict criterion is warranted so that it can handle all corner cases.

\section{Infinite eigenvalues}\label{sec:infiniteeigenvalues}

In this section, we will drop the assumption that all the eigenvalues are finite. This means that instead of an eigenvalue being a scalar, it is now a pair of scalars $(\alpha, \beta) \ne (0,0)$ so that $\det (\beta H - \alpha T) = 0$. Note that if $\beta \ne 0$, we can recover the scalar representation of the eigenvalues as the fraction $\lambda = \frac{\alpha}{\beta}$. If $\beta = 0$ we say that the eigenvalue is infinite.

In subsection~\ref{subsection:numericallyinfinitedefinition}, we will define what it means for a pencil to have a numerically infinite eigenvalue. In subsection~\ref{subsection:heuristics}, we will evaluate some commonly used criteria for infinite eigenvalues using this definition.

\subsection{Numerically infinite eigenvalues}\label{subsection:numericallyinfinitedefinition}

An eigenvalue is infinite if $\beta = 0$. This implies that the pencil has an infinite eigenvalue if and only if $T$ is singular. In practice, it is often the case that a matrix is very close to, but not exactly singular. We call such a matrix numerically singular. For our purposes, we will consider a matrix $A$ numerically singular if $\kappa_2(A) \ge u^{-1}$. This criterion is commonly used in practice. It is satisfied if and only if there exists a perturbation $E$ so that $A + E$ is exactly singular and $\|E\|_2 \le u\|A\|_2$~\cite{demmel1997applied}, i.e. $A$ is only a small perturbation away from being singular.

In a similar vein, we will say that a pencil $(H, T)$ has a numerically infinite eigenvalue if and only if there exist matrices $E_1$ and $E_2$, with $\|E_1\|_2 \le u\|H\|_2$ and $\|E_2\|_2 \le u\|T\|_2$ so that $(H+E_1, T+E_2)$ has an infinite eigenvalue. With this definition, a pencil has a numerically infinite eigenvalue if and only if $T$ is numerically singular.

For finite eigenvalues, an essential property of our new deflation criterion is that it considers values in both $H$ and $T$ and is thus invariant to diagonal scaling. If we consider an eigenvalue infinite if and only if $T$ is numerically singular, then values in $H$ are ignored. This would seem to be a valid argument against our definition, but this subsection will argue that invariance to diagonal scaling is not required.

Consider the two pencils:
\begin{equation*}
   (H_1, T_1) =
   \left(
   \begin{bmatrix}
      1 & 0 \\
      0 & 10
   \end{bmatrix},
   \begin{bmatrix}
      1 & 0 \\
      0 & 1
   \end{bmatrix}
   \right)
\end{equation*}
and
\begin{equation*}
   (H_2, T_2) =
   \left(
   \begin{bmatrix}
      1 & 0 \\
      0 & 10^{-7}
   \end{bmatrix},
   \begin{bmatrix}
      1 & 0 \\
      0 & 10^{-8}
   \end{bmatrix}
   \right).
\end{equation*}
In exact arithmetic, these pencils can be transformed into each other through row or column scaling and they both have one eigenvalue $1$ and one eigenvalue $10$. In single precision arithmetic, $T_2$ can be considered numerically singular, whereas $T_1$ is unitary and is as far from singular as possible. According to our definition, $T_1$ has only finite eigenvalues, whereas $T_2$ has a numerically infinite eigenvalue.

It is common to perform row and column scalings to improve the conditioning of the eigenvalue problem \cite{ward1981balancing,lemonnier2006balancing,dopico2022diagonal}, so it would seem that a good criterion to identify infinite eigenvalues should be invariant to such scalings. To be invariant to scaling, we need to consider the ratio of $h_{i,i}$ and $t_{i,i}$ and for the pencil $(H_2, T_2)$ this seems like a good idea. However, consider the following pencils:
\begin{equation*}
   (H_3, T_3) =
   \left(
   \begin{bmatrix}
      1 & 0 \\
      0 & 10^8
   \end{bmatrix},
   \begin{bmatrix}
      1 & 0 \\
      0 & 1
   \end{bmatrix}
   \right)
\end{equation*}
and
\begin{equation*}
   (H_4, T_4) =
   \left(
   \begin{bmatrix}
      1 & 0 \\
      0 & 1
   \end{bmatrix},
   \begin{bmatrix}
      1 & 0 \\
      0 & 10^{-8}
   \end{bmatrix}
   \right).
\end{equation*}
Just like before, these pencils can be transformed into each other through diagonal similarity transformation. They both have eigenvalues 1 and $10^8$. In single precision, $(H_3, T_3)$ has finite eigenvalues and $(H_4, T_4)$ has an infinite eigenvalue. In the first example, it may seem clear that the pencil is badly scaled, but here it is more difficult. If we decide that $(H_3, T_3)$ is the `correct' one, then any pencil with infinite eigenvalues should be converted to a pencil with large eigenvalues so long as $t_{i,i}$ is not exactly zero. If we decide that $(H_4, T_4)$ is the `correct' one, then any pencil with large eigenvalues should be converted to a pencil with infinite eigenvalues. In other words, there is no way to distinguish a badly scaled pencil from a pencil with infinite eigenvalues without additional information.

\subsection{Heuristics}\label{subsection:heuristics}

We have defined that a pencil has an infinite eigenvalue if and only if $T$ is numerically singular. However, exactly determining the condition number of $T$ in each iteration would be prohibitively expensive. The normwise criterion:
\begin{equation}
   |t_{i,i}| \le u \|T\|_F
\end{equation}
and the elementwise criterion:
\begin{equation}
   |t_{i,i}| \le u (t_{i,i-1} + t_{i+1,i}).
\end{equation}
can be interpreted as cheap heuristics to determine whether the upper triangular $T$ is numerically singular. We will study two types of errors the heuristics might make: they could identify a finite eigenvalue as numerically infinite or they could identify a numerically infinite eigenvalue as finite.

A lower bound on the condition number of a triangular matrix is: \cite{higham1987survey}
\begin{equation}\label{eq:2normconditionbound}
   \kappa_2(T) := \|T\|_2 \|T^{-1}\|_2 \ge \|T\|_2 (\min |t_{i,i}|)^{-1}.
\end{equation}
The normwise criterion does not use the 2-norm, but the Frobenius norm (because it is cheaper to compute). The Frobenius norm and the 2-norm of any matrix $A$ are related via:
\begin{equation}\label{eq:relation2andfrobnorm}
   \|A\|_2 \le \|A\|_F \le \sqrt{n}\|A\|_2.
\end{equation}
Combining Equations \eqref{eq:2normconditionbound} and \eqref{eq:relation2andfrobnorm}, we get that:
\begin{equation}
   \begin{split}
      &|t_{i,i}| \le u \|T\|_F\\
      &|t_{i,i}| \le \sqrt{n} u \|T\|_2\\
      &\Rightarrow \kappa_2(T) \ge u^{-1}n^{-\frac{1}{2}}.\\
   \end{split}
\end{equation}
So even if the normwise criterion misclassifies a finite eigenvalue as infinite, it will at least be close to numerically infinite by a factor $n^{-\frac{1}{2}}$. However, the normwise criterion may misclassify a numerically infinite eigenvalue as finite.

The elementwise criterion, which is stricter than the normwise criterion, is then also guaranteed to only identify eigenvalues that are at least close to being infinite. Because it is stricter, it is also more likely to falsely identify numerically infinite eigenvalues as finite. In the standard eigenvalue problem, an elementwise criterion was relevant to accurately determine small eigenvalues in graded matrices. However, graded matrices are often numerically singular and therefore (by our definition) lead to infinite eigenvalues so that argument is not valid here.

We now present some experiments that test both types of errors the criteria can make. They will show that (at least for these tests) neither of the criteria falsely identifies a finite eigenvalue as infinite, but sometimes they fail to identify numerically infinite eigenvalues. To perform the tests, the pencils are always reduced to Hessenberg-triangular form and then the QZ algorithm (a lightly edited version of XHGEQZ as before) is run in double precision.

First, we consider a pencil $(A,B)$ where the sparsity pattern of $B$ is such that it has a number of infinite eigenvalues.
\begin{equation}
   B = \begin{bmatrix}
      B_1 & 0\\
      0 & B_2,
   \end{bmatrix}
\end{equation}
with $B_1$ an $m_1 \times m_2$ matrix with randomly drawn entries and $B_2$ an $m_2 \times m_1$ matrix with randomly drawn entries. $A$ is a full matrix with randomly drawn entries. The pencil has $|m_2 - m_1|$ infinite eigenvalues. We test 1000 of these pencils, with $m_1 = 22$ and $m_2 = 28$, so we expect to find 6 infinite eigenvalues. Using the normwise criterion, we find always find 6 infinite eigenvalues. Using the elementwise criterion, we only find 4.41 infinite eigenvalues on average.

Next, we generate pencils using the procedure described in the previous section, with $\alpha_i = 1$, $\beta_i = 10^{-16i/n}$, and $\kappa = 1$. In double precision for matrices of size 50, the last two singular values of $B$ are $1.0\cdot 10^{-16}$ and $2.089\cdot 10^{-16}$, the first singular value is $1$ and $u = 2.2204\cdot 10^{-16}$, so we expect to find 2 numerically infinite eigenvalues. The numerically infinite eigenvalues are very close to being finite, and several finite eigenvalues are very close to being numerically infinite so this is a good stress test of the criterion. Just like before, we test 1000 of these. On average, we find 1.94 and 0 infinite eigenvalues using the normwise and elementwise criteria respectively.

In both experiments, the normwise criterion performed better. Combining the theoretical arguments and the experimental results, we conclude that the normwise criterion is a better estimate of our definition.

\subsection{Overriding the default behavior}


We have argued in Subsection~\ref{subsection:numericallyinfinitedefinition} that we should consider a pencil $(A, B)$ to have a numerically infinite eigenvalue if and only if $B$ is numerically singular. This definition has some flaws, but we have argued that a better definition is not possible without taking extra information into account. An example where such extra information is available can be found in the work of Van Barel and Tisseur~\cite{van2018polynomial, tisseur2021min}. They present a scaling method for pencils that arise when looking for roots of polynomials. This scaling significantly improves the accuracy of the computed eigenvalues. Unfortunately, it is also prone to introducing numerically infinite eigenvalues. Because of the way the pencil is constructed, they know all of these eigenvalues should be considered finite. They accomplish this by using an extremely strict deflation criterion:
\begin{equation*}
   |t_{i,i}| < u_s,
\end{equation*}
where $u_s$ is the smallest positive floating point number. This essentially means that any eigenvalue pair with $\beta \ne 0$ is considered finite.

In order to accommodate both normal users who want a reasonable default and expert users who can be certain that no infinite eigenvalues are present, we propose using an optional flag. By default, the QZ algorithm would use the normwise criterion, but by specifying the flag, the user can switch to the extra strict criterion by Van Barel and Tisseur~\cite{van2018polynomial, tisseur2021min}.

\section{Conclusion}

We derived a new deflation criterion for finite eigenvalues in the QZ algorithm. It has a better mathematical foundation and both theoretical arguments and numerical experiments indicate that it is at least as, and sometimes much more accurate than other commonly used criteria. We also evaluated commonly used criteria for infinite eigenvalues and concluded that the normwise criterion performs well as a reasonable default, but that expert users should be able to switch to an extra strict criterion if needed.


\section*{Acknowledgement}
We are grateful to the Mathworks team for reporting the issues related to deflations of infinite eigenvalues in LAPACK version 3.10 which formed the inspiration for this paper.

\subsection*{Conflicts of interest}
This study does not have any conflicts to disclose.

\subsection*{Author contributions}

J. Langou encouraged T. Steel to investigate the deflation criterion in the QZ algorithm. T. Steel worked out the technical details and performed the numerical experiments. R. Vandebril contributed to the theoretical analysis. R. Vandebril and J. Langou supervised the project. T. Steel wrote the manuscript with input and feedback from all authors.

\subsection*{Financial disclosure}

The research of T. Steel and R. Vandebril was partially supported by the Research Council KU Leuven, project C16/21/002 (Manifactor: Factor Analysis for Maps into Manifolds). The research of J. Langou was partially supported by NSF award \#2004850.

\nocite{*}
\bibliography{references}%



\end{document}